\documentclass{amsart}
\def\ds{\displaystyle}
\def\N{{\mathbb N}}
\newcommand{\R}{{\mathbb R}}
\newcommand{\D}{{\mathbb D}}
\newcommand{\SF}{{\mathbb S}}

\newcommand{\loc}{{\rm loc}}
\newcommand{\ess}{\mathop{\rm esssup}\limits}
\newcommand{\es}{\mathop{\rm essinf}\limits}

\def\XXint#1#2#3{{\setbox0=\hbox{$#1{#2#3}{\int}$ }
\vcenter{\hbox{$#2#3$ }}\kern-.6\wd0}}

\def\wtl{\widetilde}
\def\ol{\overline}%%%%%%%
\def\B{{\mathcal B}}
\def\C{{\mathcal C}}
\def\E{{\mathcal  E}}

\def\K{{\mathcal K}}
\def\P{{\mathcal P}}
\def\SS{{\mathcal S}}

\def\I{{\mathcal I}}
\def\T{{\mathcal T}}
\def\KK{{\mathfrak K}}
\def\CC{{\mathfrak C}}
\def\II{{\mathfrak I}}

\def\sm{\setminus}%

\def\ba{{\mathbf a}}
\newcommand{\supp}{{\rm supp\,}}

\newtheorem{thm}{Theorem}[section]

\newtheorem{lem}[thm]{Lemma}
\newtheorem{cor}[thm]{Corollary}

\newtheorem{defn}[thm]{Definition}

\newtheorem{example}[thm]{Example}

\begin{document}

\title[Parabolic equations in generalized Morrey spaces] {Generalized Morrey regularity for parabolic equations with discontinuity data}

\author[V.S. Guliyev, L.G. Softova]{Vagif S. Guliyev${}^1$, \quad  Lubomira G. Softova${}^2$}

\address{${}^1$Ahi Evran University\\
 Department of Mathematics\\
 Kirsehir, Turkey\\
Institute of Mathematics and Mechanics of NAS of Azerbaijan\\
 Baku}
\email{vagif\@@guliyev.com}

\address{${}^2$Department of Civil Engineering\\
  Second University  of Naples\\
 Via Roma 29\\
 81031 Aversa\\
 Italy}
\email{luba.softova\@@unina2.it}



\maketitle

\begin{abstract}
We obtain  continuity in generalized parabolic Morrey spaces  of sublinear integrals
generated by the parabolic Calder\'{o}n-Zygmund operator and its commutator with  $VMO$ functions. 
The obtained estimates are used to study global regularity of the solutions
of the Cauchy-Dirichlet   problem for linear  uniformly parabolic equations 
  with discontinuous coefficients. 
 \end{abstract}

\section{Introduction}

The classical Morrey spaces $L_{p,\lambda}$ are originally introduced  in \cite{Mo} in order
  to prove local H\"older continuity  of solutions to certain systems of  partial
differential equations (PDE's).
A real valued function $f$ is said to belong to the Morrey space $L_{p,\lambda}$ with $p\in[1,\infty),$ $\lambda\in(0,n)$  provided the following norm is finite
$$
\|f\|_{L_{p,\lambda}(\R^n)}=\left( \sup_{(x,r)\in\R^n\times\R_+ }\frac1{r^\lambda}
\int_{\B_r(x)}|f(y)|^p\,dy   \right)^{1/p}\,.
$$
The main result connected with  these spaces is the following celebrated lemma: let $|Df|\in L_{p,\lambda}$ even locally, with $\lambda<p,$ then $u$ is H\"older continuous of exponent $\alpha =1-\frac{\lambda}{p}.$
This result has found many applications in the study the regularity of the strong solutions to elliptic and parabolic PDE's and systems.
 In \cite{ChFra} Chiarenza and Frasca showed
 boundedness of the Hardy-Littlewood maximal operator in $L_{p,\lambda}(\R^n)$ that allows them to prove continuity in that spaces of some   classical
 integral operators.  These  operators  appear in the representation formulas of the solutions of linear  PDE's  and systems. Thus the results in \cite{ChFra}  permit to study the  regularity  of the solutions of these operators  in 
  $L_{p,\lambda}$   (see \cite{PS1, Sf1}).  In \cite{Mi} Mizuhara extends the concept of Morrey of integral average over a ball with a certain growth, taking a weight   function
$\omega(x,r):\R^{n+1}\times\R_+\to\R_+$ instead of $r^\lambda.$ Thus he put
the beginning of the study of the   generalized Morrey spaces $L_{p,\omega}$ under various conditions on the weight function.  In    \cite{Na}  Nakai extended  the results of \cite{ChFra}
in $L_{p,\omega}$ imposing the following conditions  on the weight
 $$
 \int_r^\infty \frac{\omega(x,s)}{s^{n+1}} ds  \leq  C \, \frac{\omega(x,r)}{r^n}, \qquad
C_1\leq \frac{ \omega(x,s)}{\omega(x,r)} \leq C_2  \  \  r \leq  s \leq 2r,
 $$
where the constants do  not depend on $s$, $r$ and $x.$
In \cite{Sf2, Sf3, Sf4} global $L_{p,\omega}$-regularity of solutions to elliptic and parabolic boundary value problems is obtained using explicit representation formula.

Other generalizations of the Morrey spaces are considered in \cite{AGM, G1, G2, GAKS} where the continuity of sublinear operators  generated by various  classical integral operators as the Calder\'on-Zygmund, Riesz and others is proved.    In \cite{GS} we have applied these results to the study of regularity of solutions to the Dirichlet problem for  linear uniformly elliptic equations.

In the present work we obtain global regularity of the solutions of the Cauchy-Dirichlet problem for parabolic non-divergence equations   with $VMO$ coefficients in $M_{p,\varphi}$.   This problem has been studied in the framework of the Morrey spaces in \cite{PRS} and in the weighted Lebesgue spaces in \cite{GM}. Here we extend these results in $M_{p,\varphi}.$ 
 For this goal we study continuity in  
$M_{p,\varphi}$ of sublinear operators generated by the Calder\'on-Zygmund  integrals with parabolic kernels and their commutators with  $BMO$ functions (Section~\ref{sec3}).   The last ones enter in the interior representation formula of the  derivatives $D_{ij}u$ of the solution of \eqref{CDP}.   
In Section~\ref{sec4}  we establish continuity  for sublinear integrals
generated by  nonsingular integral operators  and commutators. These integrals  enter in the boundary representation formula for  $D_{ij}u.$ The global a priori estimate for $u$ is obtained in Section~\ref{sec6}.

Throughout this paper the following notations will be used:
\begin{itemize}
\item $  x=(x',t), y=(y',\tau)\in \R^{n+1}=\R^n\times \R, \  \R^{n+1}_+=\R^n\times\R_+;$
\item $ x=(x'',x_n,t)\in \D^{n+1}_+=\R^{n-1}\times\R_+\times\R_+, \   \D^{n+1}_-=\R^{n-1}\times\R_-\times\R_+;$
\item $|\cdot|$ is the Euclidean metric, $  |x|=\left(\sum_{i=1}^n x_i^2 +t^2 \right)^{1/2};$
\item $D_iu=\partial u/ \partial x_i, \  Du=(D_1u,\ldots,D_nu), \  u_t=\partial u/\partial t;$
\item $  D_{ij}u= \partial^2 u/\partial x_i\partial x_j, \ D^2u=\{D_{ij}u\}_{ij=1}^n $  means the Hessian matrix of  $u;$
\item $\B_r(x')=\{y'\in \R^n:\  |x'-y'|<r\},  $ $|\B_r|=C r^n;$
\item $\I_r(x) =\{y\in \R^{n+1}: \   |x'-y'|<r, |t-\tau|<r^2 \},$ $|\I_r|=Cr^{n+2};$
\item ${\SF}^{n}$  is the unit sphere in $ {\R}^{n+1};$ 
\item for any $ f\in L_p(A),$ $A\subset \R^{n+1}$ we  write
$$
\|f\|_{p,A}\equiv \|f\|_{L_p(A)}= \left( \int_A |f(y)|^p dy  \right)^{1/p}.
$$
\item  The standard summation convention on repeated upper and lower indexes is  adopted.
\item  The letter $C$ is used for various positive  constants and may change from one occurrence to another.
\end{itemize}

\section{Definitions and statement of the problem}\label{sec2}

In the following,  besides   the standard parabolic  metric
$\varrho(x)=\max(|x'|,|t|^{1/2})$
we  use  the equivalent  one
 $ \rho(x)=\left(\frac{|x'|^2+\sqrt{|x'|^4+4t^2}}{2}\right)^{1/2} $ introduced by Fabes and Rivi\'ere in \cite{FR}.
The induced by it topology consists of  ellipsoids
$$
{\E}_r(x)=\left\{y\in \R^{n+1}:\   \frac{|x'-y'|^2}{r^2}
+\frac{|t-\tau|^2}{r^4} <1\right\}, \  |\E_r|=C r^{n+2}, \   \E_1(x)\equiv \B_1(x)\,.
$$
It is easy to see that the metrics $\rho(\cdot)$ and $\varrho(\cdot)$ are equivalent. Infact
 for each $\E_r$ there exist parabolic cylinders
$\underline \I$ and $\overline \I$ with measure comparable to $r^{n+2}$ such that ${\underline \I}\subset \E_r\subset {\overline \I}.$ In what follows
 all estimate obtained over ellipsoids hold true also over parabolic cylinders and we shall use this property without explicit references.

Let $\Omega\subset\R^n$ be a bounded $C^{1,1}$-domain and $Q=\Omega\times(0,T),$ $T>0$ be a cylinder in $\R^{n+1}_+.$
We   give the definitions of the functional spaces   which we are going to use.
\begin{defn}\label{VMO}
Let $a\in L_1^{\rm loc}(\R^{n+1})$ and  $ a_{\E_r}=|\E_r|^{-1}\int_{\E_r} a(y)dy$ be the mean integral of $a.$ Denote
$$
\eta_a(R)=\sup_{r\leq R}\frac{1}{|\E_r|}\int_{\E_r}|f(y)-f_{\E_r}| dy \quad \text{ for every } R>0
$$
where $\E_r$ ranges over all ellipsoids in $\R^{n+1}.$ 
We say that
\begin{enumerate}
\item[$\bullet$]
$a\in BMO$  ({\it bounded mean oscillation},  \cite{JN})  provided the following is finite
$$
\|a\|_{\ast}=\sup_{R>0}\eta_a(R).
$$
The quantity $\|\cdot\|_\ast$ is a norm in $BMO$ modulo constant function under which $BMO$ is a Banach space.
\item[$\bullet$] $a\in VMO$   ({\it vanishing mean oscillation},  \cite{Sar})  if $a\in BMO$ and
$$
\lim_{R\to 0}\eta_a(R)=0.
 $$
 The quantity $\eta_a(R)$ is called $VMO$-modulus of $a.$
\end{enumerate}

For any bounded cylinder  $Q$ we define
 $BMO(Q)$ and $VMO(Q)$ taking  $a\in L_1(Q)$ and  $Q_r $  instead of $\E_r$ in the definition above.
\end{defn}

According to \cite{A, Jones}, having a function $a\in BMO(Q)$ or  $VMO(Q)$ it is possible to extend it in the whole ${\R}^{n+1}$ preserving its
 $BMO$-norm or $VMO$-modulus, respectively. In the following we use this property without explicit references.
Any bounded uniformly continuous (BUC) function  $f$ with modulus of continuity $\omega_f(R)$ belongs to $VMO$  with $\eta_f(R)= \omega_f(R).$    Besides that, $BMO$ and $VMO$ contain also discontinuous functions  and the following example shows the inclusion
$W_{1,n+2}(\R^{n+1})\subset VMO\subset BMO.$
\begin{example}
$f_\alpha(x)= |\log\rho(x)|^\alpha \in VMO$ for any $\alpha\in(0,1);$\\
  $f_\alpha\in W_{1,n+2}({\R}^{n+1})$
for $\alpha\in(0,1-1/(n+2));$\\
 $f_\alpha\notin W_{1,n+2}({\R}^{n+1})$ for $\alpha\in[1-1/(n+2), 1);$\\
$f(x)=|\log\rho(x)|\in BMO\setminus VMO;$\\
$\sin f_\alpha(x)\in VMO\cap L_\infty(\R^{n+1}).$
\end{example}

\begin{defn} \label{def1}
Let $\varphi:\R^{n+1}\times\R_+\to\R_+$ be a measurable function and $ p\in[1,\infty).$ The generalized parabolic Morrey space $M_{p,\varphi}(\R^{n+1})$ consists of all functions $f\in L_p^{\loc}(\R^{n+1})$
such that
$$
\|f\|_{p,\varphi;\R^{n+1}} =\sup_{(x,r)\in\R^{n+1}\times\R_+}
\varphi(x,r)^{-1} \left(r^{-(n+2)} \, \int_{\E_r(x)}|f(y)|^p dy\right)^{1/p}<\infty.
$$
The space  $M_{p,\varphi}(Q)$  consists of  $L_p(Q)$ functions provided the following norm is finite
$$
\|f\|_{p,\varphi;Q} =\sup_{(x,r)\in Q\times\R_+}
\varphi(x,r)^{-1} \left(r^{-(n+2)} \, \int_{Q_r(x)}|f(y)|^p dy\right)^{1/p}
$$
where $Q_r(x)=Q\cap \I_r(x).$
The generalized weak parabolic Morrey space $WM_{1,\varphi}(\R^{n+1})$ consists of all measurable functions such that
$$
\|f\|_{WM_{1,\varphi}(\R^{n+1})}=\sup_{(x,r)\in\R^{n+1}\times\R_+}
\varphi(x,r)^{-1}r^{-n-2}
\| f\|_{WL_1(\E_r(x))}
$$
where $WL_1$ denotes the weak $L_1$ space.

The generalized  Sobolev-Morrey  space $W^{2,1}_{p,\varphi}(Q),$ $p\in[1,\infty)$  consist of all Sobolev  functions $u\in W^{2,1}_{p}(Q)$
with distributional derivatives $D^l_tD^s_xu\in M_{p,\varphi}(Q),$  $0\leq 2l+ |s|\leq 2$ endowed by the norm
$$
\|u\|_{W^{2,1}_{p,\varphi}(Q)}=\|u_t\|_{p,\varphi;Q}+ \sum_{ |s|\leq 2}\|D^s u\|_{p,\varphi;Q}.
$$
$$
 \overset{\circ}{W}{}^{2,1}_{p,\varphi}(Q)=\big\{ u\in W^{2,1}_{p,\varphi}(Q):\  u(x)=0, x\in
 \partial Q\big\}, \   
\|u\|_{\overset{\circ}{W}{}^{2,1}_{p,\varphi}(Q)}=\|u\|_{W^{2,1}_{p,\varphi}(Q)}
$$
where $\partial Q$ means the parabolic boundary $\Omega\cup(\partial\Omega\times(0,T)).$
\end{defn}

We consider the Cauchy-Dirichlet problem for linear parabolic equation
\begin{equation}\label{CDP}
\begin{cases}
u_t- a^{ij}(x)D_{ij}u(x)=f(x) \   \text{ a.a. }  x\in Q,\quad
  u\in\overset{\circ}{W}{}^{2,1}_{p,\varphi}(Q)
\end{cases}
\end{equation}
where  the coefficient matrix $\ba(x)=\{a^{ij}(x)\}_{i,j=1}^n$ satisfies
\begin{equation}\label{coef}
\begin{cases}
\exists \   \Lambda>0:\
\Lambda^{-1}|\xi|^2\leq a^{ij}(x)\xi_i\xi_j \leq \Lambda |\xi|^2\   \text{ for a.a. } x\in Q,\
\forall \xi\in \R^n\\
a^{ij}(x)=a^{ji}(x) \text{ that implies }  a^{ij}\in L_\infty(Q).
\end{cases}
\end{equation}
\begin{thm}{\bf (Main result)}\label{main}
Let $\ba\in VMO(Q)$ satisfy \eqref{coef}  and for each $p\in(1,\infty),$ $u\in \overset{\circ}{W}{}^{2,1}_{p}(Q)$   be  a strong solution of \eqref{CDP}. If $f\in M_{p,\varphi}(Q)$ with $\varphi(x,r)$ being measurable positive function  satisfying
\begin{equation}\label{weight}
\int_{r}^{\infty} \Big(1+\ln \frac{s}{r}\Big)
\frac{\es_{s<\zeta<\infty} \varphi(x,\zeta) \zeta^{\frac{{n+2}}{p}}}{s^{\frac{{n+2}}{p}+1}}\, ds
\le C \,\varphi(x,r),\quad (x,r)\in Q\times\R_+
\end{equation}
then $u\in \overset{\circ}{W}{}^{2,1}_{p,\varphi}(Q)$ and
\begin{equation}\label{apriori}
\|u\|_{\overset{\circ}{W}{}^{2,1}_{p,\varphi}(Q)}\leq C \|f\|_{p,\varphi;Q}
\end{equation}
with $C=C(n,p,\Lambda, \partial\Omega, T, \|\ba\|_{\infty;Q}, \eta_a).$
\end{thm}

\section{ Sublinear  operators generated by parabolic singular integrals  in generalized  Morrey spaces}\label{sec3}

Let $f\in L_1(\R^{n+1})$ be a function with a compact support and $a\in BMO.$ For any  $x\notin \supp f$ define the sublinear operators $T$ and
 $T_a$ such that
\begin{align}\label{subl}
|Tf(x)|&\le C \int_{\R^{n+1}}
\frac{|f(y)|}{\rho(x-y)^{n+2}} \,dy\\
\label{sublcomm}
|T_{a}f(x)|&\le C \int_{\R^{n+1}} |a(x)-a(y)| \,
\frac{|f(y)|}{\rho(x-y)^{n+2}} \, dy.
\end{align}
Suppose in addition that the both operators are bounded in $L_p(\R^{n+1})$ satisfying the estimates
\begin{equation}\label{bounds}
\|Tf\|_{p;\R^{n+1}}\leq C\|f\|_{p;\R^{n+1}}, \quad \|T_af\|_{p;\R^{n+1}}\leq C\|a\|_{\ast}\|f\|_{p;\R^{n+1}}
\end{equation}
with constants independent of $a$ and $f.$
The following known  result concerns   the   Hardy operator
$
Hg(r)=\frac{1}{r}\int_0^r g(s)ds,
$ $r>0.$
\begin{thm}\label{thmHardy}{\rm(\cite{CarPickSorStep})}
The inequality
\begin{equation}\label{hardy}
\ess_{r>0}w(r) Hg(r) \leq A \ess_{r>0}v(r)g(r)
\end{equation}
holds for all  non-increasing functions  $g:\R_+\to \R_+$ if and
only if
\begin{equation}\label{eqA}
A= C \sup_{r>0}\frac{w(r)}{r}\int_0^r
\frac{ds}{\ess_{0<\zeta<s}v(\zeta)}<\infty.
\end{equation}
\end{thm}
\begin{lem}\label{lem1}
Let   $f\in L_p^{\rm loc}(\R^{n+1}),$  $p\in[1,\infty)$ be such that 
\begin{equation}\label{eq2}
\int_r^{\infty} s^{-\frac{n+2}{p}-1} \|f\|_{p;\E_s(x_0)}ds<\infty \quad \forall \      
   (x_0,r)\in\R^{n+1}\times\R_+
\end{equation}
and $T$    be a  sublinear operator satisfying  \eqref{subl}.

{\rm (i) }   If $p>1$ and $T$ bounded on $L_p(\R^{n+1})$ then
\begin{equation}\label{eq3}
 \|T f\|_{p;\E_r(x_0)} \le C \,
r^{\frac{n+2}{p}} \int_{2r}^{\infty} s^{-\frac{n+2}{p}-1} \|f\|_{p;\E_s(x_0)}ds.
\end{equation}

{\rm (ii) } If $p=1$ and $T$ bounded from $L_1(\R^{n+1})$ on $WL_1(\R^{n+1})$ then
\begin{equation}\label{eq3a}
\|Tf\|_{WL_1(\E_r(x_0))}\leq Cr^{n+2}\int_{2r}^\infty s^{-n-3}\|f\|_{1,\E_s(x_0)}\,ds
\end{equation}
where the constants are independent of $r,$ $x_0$ and $f.$
\end{lem}
\proof{\rm (i) } Fix  a point $x_0\in \R^{n+1}$ and  consider an ellipsoid $\E_r(x_0).$ Denote by $2\E_r(x_0)=\E_{2r}(x_0)$ and $\E_r^c(x_0)=\R^{n+1}\setminus \E_r(x_0).$
Consider the decomposition of $f$ with respect to the    ellipsoid $\E_r(x_0)$
$$
f= f\chi_{2\E_r(x_0)} + f\chi_{2\E_r^c(x_0)}=f_1+f_2.
$$
Because of the  $(p,p)$-boundedness of the operator $T$ and $f_1\in L_p(\R^{n+1})$ we have
\begin{equation*}
\|Tf_1\|_{p;\E_r(x_0)}\leq \|Tf_1\|_{p;\R^{n+1}}\leq
C\|f_1\|_{p;\R^{n+1}}=C\|f\|_{p;2\E_r(x_0)}.
\end{equation*}
It is easy to see that for arbitrary points  $x\in \E_r(x_0)$ and  $y\in 2\E_r^c(x_0)$  it holds
\begin{equation}\label{xy}
\frac{1}{2}\rho(x_0-y)\le\rho(x-y)\le \frac{3}{2}\rho(x_0-y).
\end{equation}
Applying \eqref{subl}, \eqref{xy}, the Fubini theorem  and the H\"older inequality  to $Tf_2$ we get
\begin{align*}
|Tf_2(x)|\leq &   C \int_{2\E_r^c(x_0)}\frac{|f(y)|}{\rho(x_0-y)^{{n+2}}}dy\leq C\int_{2\E_r^c(x_0)}|f(y)|\left(\int_{\rho(x_0-y)}^\infty \frac{ds}{s^{{n+3}}}\right)dy\\[6pt]
\leq & C\int_{2r}^\infty \left(\int_{2r\leq \rho(x_0-y)<s}|f(y)|dy  \right)\frac{ds}{s^{{n+3}}}\\[6pt]
\leq & C \int_{2r}^\infty \left(\int_{\E_s(x_0)}|f(y)|dy  \right)\frac{ds}{s^{{n+3}}}\leq C \int_{2r}^\infty \|f\|_{p;\E_s(x_0)}\frac{ds}{s^{\frac{n+2}{p}+1}}\, .
\end{align*}
Direct calculations give
\begin{equation} \label{eq4}
\|Tf_2\|_{p,\E_r(x_0)}\leq C r^{\frac{{n+2}}{p}} \int_{2r}^{\infty}\|f\|_{p;\E_s(x_0)}\frac{ds}{s^{\frac{{n+2}}{p}+1}}
\end{equation}
which holds for all  $p\in[1,\infty).$ Thus
\begin{equation}\label{Tf}
\|Tf\|_{p;\E_r(x_0)}\leq  C\left(  \|f\|_{p;2\E_r(x_0)}+
r^{\frac{{n+2}}{p}}\int_{2r}^{\infty}\|f\|_{p;\E_s(x_0)}\frac{ds}{s^{\frac{{n+2}}{p}+1}}\right).
\end{equation}
On the other hand
\begin{equation}  \label{eq5}
\|f\|_{p,2\E_r(x_0)}  \leq C r^{\frac{{n+2}}{p}}
\int_{2r}^{\infty}\|f\|_{p;\E_s(x_0)}\frac{ds}{s^{\frac{{n+2}}{p}+1}}
\end{equation}
which unified with \eqref{Tf} gives \eqref{eq3}.

{\rm (ii) } Let $f\in L_1(\R^{n+1}),$ the weak $(1,1)$-boundedness of $T$ implies
\begin{align*}
\|Tf_1\|_{WL_1(\E_r(x_0))}&\leq \|Tf_1\|_{WL_1(\R^{n+1})}\\
&\leq C\|f_1\|_{1,\R^{n+1}}= C \|f\|_{1,2\E_r(x_0)}\\
&\leq C r^{n+2}\int_{2r}^{+\infty}\|f\|_{1,\E_s(x_0)}\frac{ds}{s^{n+3}}
\end{align*}
that unified with \eqref{eq4} gives \eqref{eq3a}.
\endproof
\begin{thm}\label{thm1}
Let    $p\in[1,  \infty),$   $\varphi(x,r)$ be a measurable positive function satisfying
\begin{equation}\label{eq6}
\int_r^{\infty} \frac{\es_{s<\zeta<\infty}
\varphi(x,\zeta) \zeta^{\frac{{n+2}}{p}}}{s^{\frac{{n+2}}{p}+1}}\, ds \le
 C \,\varphi(x,r)\quad \forall \  (x,r)\in \R^{n+1}\times\R_+
\end{equation}
 and $T$ be sublinear operator satisfying  \eqref{subl}.

{\rm (i)} If $p>1$ and $T$ bounded on $L_p(\R^{n+1})$ than $T$ is bounded on $M_{p,\varphi}(\R^{n+1})$ and
\begin{equation}\label{normTf}
\|Tf\|_{p,\varphi;\R^{n+1}}\leq C\|f\|_{p,\varphi;\R^{n+1}}\,.
\end{equation}

{\rm (ii)} If $p=1$ and $T$ bounded from $L_1(\R^{n+1})$ to $WL_1(\R^{n+1})$ than it is bounded from
 $M_{1,\varphi}(\R^{n+1})$ to $WM_{1,\varphi}(\R^{n+1})  $ and
\begin{equation}\label{normTf1}
\|Tf\|_{WM_{1,\varphi}(\R^{n+1})}\leq C \|f\|_{1,\varphi;\R^{n+1}}
\end{equation}
with  constants  independent on $f$.
\end{thm}
\proof  {\rm (i) }
By  Lemma~\ref{lem1} we have
\begin{align*}
\|Tf\|_{p,\varphi;\R^{n+1}} & \leq C
\sup_{(x,r)\in\R^{n+1}\times\R_+}\varphi(x,r)^{-1}
\int_r^{\infty}\|f\|_{p;\E_s(x)}\,\frac{ds}{s^{\frac{{n+2}}{p}+1}}\\
& =C \sup_{(x,r)\in\R^{n+1}\times\R_+}\varphi(x,r)^{-1}
\int_0^{r^{-(n+2)/p}}\|f\|_{p;\E_{s^{-p/(n+2)}}(x)}\,ds\\
&= C \sup_{(x,r)\in\R^{n+1}\times\R_+}\varphi(x,r^{-p/(n+2)})^{-1}
\int_0^r \|f\|_{p;\E_{s^{-p/(n+2)}}(x)}\, ds.
\end{align*}
Applying the Theorem~\ref{thmHardy} with
\begin{align*}
& w(r)=v(r)=r\varphi(x,r^{-p/(n+2)})^{-1}, \quad   g(r)=\|f\|_{p;\E_{r^{-p/(n+2)}}(x)},\\
& Hg(r)=r^{-1}\int_0^r \|f\|_{p;\E_{s^{-p/(n+2)}}(x)}\,ds,
\end{align*}
 where  the condition  \eqref{eqA} is equivalent to   \eqref{eq6},   we get \eqref{normTf}.

{\rm (ii) } Making use of \eqref{eq3a} and \eqref{hardy} we get
\begin{align*}
&\|Tf\|_{WM_{1,\varphi}(\R^{n+1})}\leq C \sup_{(x_0,r)\in\R^{n+1}\times\R_+} \varphi(x_0,r)^{-1} 
\int_r^\infty\|f\|_{1,\E_s(x_0)} \frac{ds}{s^{n+3}}\\
&\quad = C \sup_{(x_0,r)\in\R^{n+1}\times\R_+}\varphi(x_0,r^{-\frac{1}{n+2}})^{-1} \int_0^r
 \|f\|_{1,\E_{s^{-1/(n+2)}}(x_0)}\,ds\\
&\quad  \leq  C \sup_{(x_0,r)\in\R^{n+1}\times\R_+}\varphi(x_0,r^{-\frac{1}{n+2}})^{-1}
r \|f\|_{1,\E_{r^{-1/(n+2)}}(x_0)}=C\|f\|_{1,\varphi;\R^{n+1}}\,.
\end{align*}
\endproof

 Our next step is to  show boundedness of $T_a$ in $M_{p,\varphi}(\R^{n+1}).$
For this goal we recall some  properties of the $BMO$ functions.
\begin{lem}{\rm(John-Nirenberg  lemma, \cite[Lemma~2.8]{BC})}\label{lemJN}
Let $a\in BMO$ and $p\in [1,\infty)$. Then for any  $\E_r$ there
holds
$$
\left( \frac{1}{|\E_r|}\int_{\E_r}|a(y)-a_{\E_r}|^p dy\right)^{\frac{1}{p}}
\leq C(p) \|a\|_{*} .
$$
\end{lem}

As an immediate consequence of Lemma~\ref{lemJN} we get the following property.
\begin{cor}
 Let $a\in BMO$ then for all $0<2r<s$  it  holds
\begin{equation}\label{propBMO}
\left|a_{\E_r}-a_{\E_s}\right| \le C(n)\big(1+ \ln \frac{s}{r}  \big) \|a\|_\ast\,.
\end{equation}
\end{cor}
\proof
Since $s>2r$ there exists $k\in {\N},$ $k\geq 1$ such that
$2^kr<s\leq 2^{k+1}r$ and hence $k\ln 2<\ln \frac{s}{r}\leq (k+1)\ln 2.$  By \cite[Lemma~2.9]{BC} we have
\begin{align*}
|a_{\E_s} -a_{\E_r} |&\leq |a_{2^k\E_r}-a_{\E_r}| + |a_{2^k\E_r} -a_{\E_s}|\\
&\leq C(n)k\|a\|_\ast+ \frac1{|2^k\E_r|} \int_{2^k\E_r}  |a(y)-a_{\E_s}|dy\\
&\leq C(n)\left(k \|a\|_\ast + \frac{1}{|\E_s|}\int_{\E_s}  |a(y)-a_{\E_s}|dy\right)\\
&<C(n)\big(\ln\frac{s}{r}+1\big)\|a\|_\ast\,.
\end{align*}
\endproof

To estimate the norm  of $T_a$ we shall employ the same idea which we 
used in the proof of Lemma~\ref{lem1}.
\begin{lem}\label{lem2}
Let  $a \in BMO$ and  $T_{a}$ be a bounded  operator in $L_p(\R^{n+1}),$ $p\in(1,\infty)$ satisfying
 \eqref{sublcomm} and \eqref{bounds}. Suppose that for  any
$f\in L_p^\loc(\R^{n+1})$ 
\begin{equation}
\int_r^{\infty} \Big(1+\ln \frac{s}{r}\Big)
\|f\|_{p;\E_s(x_0)}\,\frac{ds}{s^{\frac{{n+2}}{p}+1} }<\infty\quad \forall \  (x_0,r)\in\R^{n+1}\times\R_+\,.
\end{equation}
Then
\begin{equation}\label{eq1}
\|T_{a}f\|_{p;\E_r(x_0)} \leq C \|a\|_{*} \, r^{\frac{{n+2}}{p}}
\int_{2r}^{\infty}  \Big(1+\ln \frac{s}{r}\Big)  \,
\|f\|_{p;\E_s(x_0)}\,\frac{ds}{s^{\frac{{n+2}}{p}+1}}
\end{equation}
where $C$ is independent of $a$, $f$, $x_0$ and $r$.
\end{lem}
\proof
Fix a point $x_0\in \R^{n+1}$ and 
consider the decomposition
$f= f\chi_{2\E_r(x_0)}+  f\chi_{2\E_r^c(x_0)}=  f_1+f_2.$  Hence
$$
\|T_{a}f\|_{p;\E_r(x_0)}\leq \|T_{a}f_1\|_{p;\E_r(x_0)}+\|T_{a}f_2\|_{p;\E_r(x_0)}
$$
and  by \eqref{bounds} as in Lemma~\ref{lem1} we have
\begin{equation}\label{eq7}
\|T_{a}f_1\|_{p;\E_r(x_0)}\leq   C\|a\|_{*} \, \|f\|_{p;2\E_r(x_0)}.
\end{equation}
On the other hand, because of \eqref{xy} we can write
\begin{equation*}
\begin{split}
\|T_{a}f_2\|_{p;\E_r(x_0)}&\leq C \left(\int_{\E_r(x_0)}\left(\int_{2\E_r^c(x_0)}
\frac{|a(x)-a(y)||f(y)|}{\rho(x_0-y)^{{n+2}}}dy\right)^pdx\right)^{\frac{1}{p}}
\\
&\leq C \left(\int_{\E_r(x_0)}\left(\int_{2\E_r^c(x_0)}
\frac{|a(y)-a_{\E_r(x_0)}||f(y)|}{\rho(x_0-y)^{{n+2}}}dy\right)^pdx\right)^{\frac{1}{p}}
\\
&+ C\left(\int_{\E_r(x_0)}\left(\int_{2\E_r^c(x_0)}\frac{|a(x)-
a_{\E_r(x_0)}||f(y)|}{\rho(x_0-y)^{{n+2}}}dy\right)^pdx\right)^{\frac{1}{p}} \\
&=I_1+I_2.
\end{split}
\end{equation*}
Applying \eqref{sublcomm}, the Fubini theorem and the H\"older inequality as in Lemmate~\ref{lem1} and~\ref{lemJN} we get
\begin{align*}
I_1 &\leq C  r^{\frac{{n+2}}{p}}\left( \int_{2r}^{\infty}\int_{\E_s(x_0)}|a(y)-a_{\E_r(x_0)}||f(y)|dy\right) \frac{ds}{s^{{n+3}}}\\
& \leq C r^{\frac{{n+2}}{p}}\left( \int_{2r}^{\infty}\int_{\E_s(x_0)}
|a(y)-a_{\E_s(x_0)}||f(y)|dy\right)\frac{ds}{s^{{n+3}}}\\
&\quad + C r^{\frac{{n+2}}{p}} \int_{2r}^{\infty}|a_{\E_r(x_0)}-a_{\E_s(x_0)}|
\left( \int_{\E_s(x_0)} |f(y)|dy\right)\frac{ds}{s^{{n+3}}}
\\
&\leq C r^{\frac{{n+2}}{p}} \int_{2r}^{\infty}
\left(\int_{\E_s(x_0)}|a(y)-a_{\E_s(x_0)}|^{\frac{p}{p-1}}dy\right)^{\frac{p-1}{p}}
\|f\|_{p;\E_s(x_0)}\,\frac{ds}{s^{{n+3}}}\\
&\quad + C r^{\frac{{n+2}}{p}} \int_{2r}^{\infty}|a_{\E_r(x_0)}-a_{\E_s(x_0)}|
\|f\|_{p;\E_s(x_0)}\,\frac{ds}{s^{\frac{{n+2}}{p}+1}}
\\
& \leq C \|a\|_{*}\,r^{\frac{{n+2}}{p}}
\int_{2r}^{\infty}\Big(1+\ln \frac{s}{r}\Big)
\|f\|_{p;\E_s(x_0)}\,\frac{ds}{s^{\frac{{n+2}}{p}+1}}.
\end{align*}
In order to estimate $I_2$ we note that
\begin{equation*}
I_2 = \left(\int_{\E_r(x_0)}|a(x)-a_{\E_r(x_0)}|^pdx\right)^{\frac{1}{p}} \int_{2\E_r^c(x_0)}
 \frac{|f(y)|}{\rho(x_0-y)^{n+2}}dy.
\end{equation*}
By Lemma~\ref{lemJN} and  \eqref{eq4}  we get
$$
I_2\leq C\|a\|_\ast\,r^{\frac{{n+2}}{p}}\int_{2\E_r^c(x_0)}
\frac{|f(y)|}{\rho(x_0-y)^{{n+2}}}dy\leq
 C \|a\|_\ast\,r^{\frac{{n+2}}{p}}
\int_{2r}^{\infty}\|f\|_{p;\E_s(x_0)}\frac{ds}{s^{\frac{{n+2}}{p}+1}}.
$$
Summing up \eqref{eq7}, $I_1$ and $I_2$  we get
$$
\|T_{a}f\|_{p;\E_r(x_0)}\leq C \|a\|_\ast\left(\|f\|_{p;2\E_r(x_0)}+
r^{\frac{{n+2}}{p}}
\int_{2r}^{\infty} \Big(1+\ln \frac{s}{r}\Big)
\|f\|_{p;\E_s(x_0)}\frac{ds}{s^{\frac{{n+2}}{p}+1}}\right)
$$
and the statement  follows after applying \eqref{eq5}.
\endproof

\begin{thm} \label{th2}
Let $p\in(1,\infty)$ and $\varphi(x,r)$ be measurable positive  function such that
\begin{equation}\label{eq8}
\int_{r}^{\infty} \Big(1+\ln \frac{s}{r}\Big)
\frac{\es_{s<\zeta<\infty} \varphi(x,\zeta) \zeta^{\frac{{n+2}}{p}}}{s^{\frac{{n+2}}{p}+1}}\, ds
\le C \,\varphi(x,r), \quad \forall \   (x,r)\in  \R^{n+1}\times\R_+
\end{equation}
where $C$ is independent of $x$ and $r$. 
Suppose $a\in BMO$ and  $T_{a}$ be  sublinear operator satisfying  \eqref{sublcomm}.
If $T_a$ is bounded in $L_p(\R^{n+1})$, then it is   bounded in
 $M_{p,\varphi}(\R^{n+1})$ and
\begin{equation}\label{normTaf}
\|T_{a}f\|_{p,\varphi; \R^{n+1}} \leq C \|a\|_\ast \, \|f\|_{p,\varphi; \R^{n+1}}
\end{equation}
with a constant independent of $a$ and $f$.
\end{thm}

The statement of the  theorem follows by Lemma~\ref{lem2} and Theorem~\ref{thmHardy}
in the same manner as the Theorem~\ref{thm1}.
\begin{example}
The functions $\varphi(x,r)=r^{\beta-\frac{n+2}{p}}$ and $\varphi(x,r)=
r^{\beta-\frac{n+2}{p}}\log^m(e+r)$ with $0<\beta<\frac{n+2}{p}$ and  $m\geq 1$ are weight functions satisfying  the condition \eqref{eq8}.
\end{example}

\section{Sublinear  operators generated by nonsingular integrals in generalized  Morrey spaces}\label{sec4}

For any  $x\in \D^{n+1}_+$ define  $\wtl x=(x'',-x_n,t)\in \D^{n+1}_-$ and 
$x^0=(x'',0,0)\in \R^{n-1}.$ Consider the semi-ellipsoids
$\E_r^+(x^0)=\E_r(x^0)\cap \D^{n+1}_+.$
Let   $f\in L_1(\D^{n+1}_+),$  $a\in BMO(\D^{n+1}_+)$ and  $\wtl T$ and  $\wtl T_a$ be  sublinear operators  such that
\begin{align}\label{tlT}
|\wtl T f(x)|&\le C \int_{\D^{n+1}_+}
\frac{|f(y)|}{\rho(\wtl x-y)^{n+2}} \,dy\\\
\label{tlTa}
|\wtl T_{a}f(x)|&\le C \int_{\D^{n+1}_+} |a(x)-a(y)| \,
\frac{|f(y)|}{\rho(\wtl x-y)^{n+2}} dy.
\end{align}
Suppose in addition that the both operators are bounded in $L_p(\D^{n+1}_+)$ satisfying the estimates
\begin{equation}\label{tlbounds}
\|\wtl T f\|_{p;\D^{n+1}_+}\leq C\|f\|_{p;\D^{n+1}_+}, \quad \|\wtl T_a f\|_{p;\D^{n+1}_+}\leq C\|a\|_{\ast}\|f\|_{p;\D^{n+1}_+}
\end{equation}
with constants independent of $a$ and $f.$
The following  assertions can be proved  in the same manner as in $\S\ref{sec3}.$
\begin{lem}  \label{lem3.3.}
Let   $f\in L_p^{\rm loc}(\D^{n+1}_+)$,  $p\in(1,\infty)$ and for all  $(x^0,r)\in \R^{n-1}\times\R_+$
\begin{equation}\label{eq9}
\int_r^{\infty} s^{-\frac{n+2}{p}-1} \|f\|_{p;\E^+_s(x^0)}ds<\infty.
\end{equation}
 If $\wtl T$ is bounded on $L_p(\D^{n+1}_+)$  then
\begin{equation}\label{eq3.5.}
\|\wtl T f\|_{p;\E^+_r(x^0)} \le C \,
r^{\frac{n+2}{p}} \int_{2r}^{\infty} s^{-\frac{n+2}{p}-1} \|f\|_{p;\E^+_s(x^0)}ds
\end{equation}
where the constant $C$ is independent of  $r,$ $x^0,$ and $f$.
\end{lem}
\begin{thm}\label{thm3}
Let   $\varphi$ be a weight function satisfying
\eqref{eq6}  and
 $\wtl T$ be a sublinear operator satisfying \eqref{tlT} and \eqref{tlbounds}. Then
it is bounded  in $M_{p,\varphi}(\D^{n+1}_+),$ $p\in(1,\infty)$ and
\begin{equation}\label{tlnormTf}
\|\wtl T f\|_{p,\varphi;\D^{n+1}_+} \leq C \|f\|_{p,\varphi;\D^{n+1}_+}
\end{equation}
with a  constant $C$ independent of $f$.
\end{thm}
\begin{lem}\label{lem5.1.}
Let $p\in(1,\infty),$ $a \in BMO(\D^{n+1}_+)$, and  ${\wtl T}_{a}$  satisfy \eqref{tlTa}
and
 \eqref{tlbounds}. Suppose that for all
$f\in L_p^{\rm loc}(\D^{n+1}_+)$ 
\begin{equation}\label{eq10}
\int_r^{\infty} \Big(1+\ln \frac{s}{r}\Big) s^{-\frac{n+2}{p}-1} \|f\|_{p;\E^+_s(x^0)}ds <\infty \quad\forall \  (x^0,r)\in\R^{n-1}\times\R_+\,.
\end{equation}
Then
$$
\|{\wtl T}_{a} f\|_{p;\E^+_r(x^0)} \leq C \|a\|_\ast  r^{\frac{n+2}{p}}
\int_{2r}^{\infty}  \Big(1+\ln \frac{s}{r}\Big)  \|f\|_{p;\E^+_s(x^0)}\,
\frac{ds}{s^{\frac{n+2}{p}+1}}
$$
with a  constant $C$ independent of $a$,$f$, $x^0$ and $r$.
\end{lem}
\begin{thm} \label{th4}
Let $p\in(1,\infty)$, $a \in BMO(\D^{n+1}_+),$ $\varphi(x^0,r)$ be a weight function satisfying   \eqref{eq8}
and  ${\wtl T}_{a}$ be a sublinear operator  satisfying  \eqref{sublcomm} and \eqref{bounds}.
Then  ${\wtl T}_{a}$ is bounded in $M_{p,\varphi}(\D^{n+1}_+),$ and
\begin{equation}\label{eq11}
\|{\wtl T}_{a} f\|_{p,\varphi;\D^{n+1}_+} \leq C \|a\|_\ast \,
\|f\|_{p,\varphi;\D^{n+1}_+}
\end{equation}
with a constant $C$ independent of $a$ and $f$.
\end{thm}

\section{Singular and nonsingular integrals  in  generalized  Morrey  spaces}\label{sec5}

In the present section we apply the above results to Calder\'on-Zygmund type operators with parabolic kernel.
  Since these operators are sublinear and bounded  in $L_p(\R^{n+1})$ their continuity in $M_{p,\varphi}$ follows immediately.
\begin{defn}\label{CZK}
A measurable function $\K(x,\xi):\R^{n+1}\times\R^{n+1}\sm\{0\}\to \R$ is called  variable parabolic Calder\'on-Zygmund kernel if:
\begin{itemize}
\item[$i)$] $\K(x,\cdot)$ is a  parabolic
 Calder\'on-Zygmund kernel for a.a. $x\in\R^{n+1}:$
\begin{itemize}
\item[$a)$] $\K(x,\cdot)\in C^\infty(\R^{n+1}\sm\{0\}),$
\item[$b)$] $\K(x,\mu\xi)=\mu^{-n-2}\K(x,\xi)$\quad $\forall \mu>0,$
 \item[$c)$] $\ds \int_{\SF^{n}}\K(x,\xi)d\sigma_\xi=0\,,$\quad  $\ds \int_{\SF^{n}}|\K(x,\xi)|d\sigma_\xi<+\infty.$
\end{itemize}
\item[$ii)$]  $\ds \left\|D^\beta_\xi \K \right\|_{L_\infty(\R^{n+1}\times\SF^{n})}\leq M(\beta)<\infty$ for every multi-index $\beta.$
\end{itemize}
\end{defn}
Moreover
$$
|\K(x,x-y)|\leq \rho(x-y)^{-n-2}\big|\K\big(x,\frac{x-y}{\rho(x-y)}\big)\big|\leq
\frac{M}{\rho(x-y)^{n+2}}
$$
which means  that the singular integrals
\begin{equation}\label{sing}
\begin{cases}
\ds\KK f(x)= P.V.\int_{\R^{n+1}}\K(x,x-y)f(y)dy\\[8pt]
\ds \CC[a, f](x)= P.V.\int_{\R^{n+1}}\K(x,x-y)[a(y)-a(x)]f(y)dy
\end{cases}
\end{equation}
are sublinear and  bounded in $L_p(\R^{n+1})$ according to the results in  \cite{BC,FR}.
Let us note that any weight function $\varphi$ satisfying
\eqref{eq8} satisfies also \eqref{eq6} and hence the
following  holds as a simple application of the estimates proved in \S\ref{sec3}.
\begin{thm}\label{CZcont}
For any $f\in M_{p,\varphi}(\R^{n+1})$ with $(p,\varphi)$ as in Theorem~\ref{th2}
and $a\in BMO$   there exist  constants  depending on $n,p$ and the kernel such that
\begin{equation} \label{sal22}
\|\KK f\|_{p,\varphi;\R^{n+1}}
\leq C\|f\|_{p,\varphi;\R^{n+1}},\quad  \|\CC[a,f]\|_{p,\varphi;\R^{n+1}}
\leq C\|a\|_\ast\|f\|_{p,\varphi;\R^{n+1}}.
\end{equation}
\end{thm}
\begin{cor}\label{locest1}
Let $Q$ be a cylinder in $\R^{n+1}_+,$  $ f\in M_{p,\varphi}(Q),$  $a\in BMO(Q)$ and
 $\K(x,\xi):\, Q\times \R^{n+1}_+\sm \{0\}\to \R.$  Then the operators \eqref{sing} are bounded in
 $M_{p,\varphi}(Q)$ and
\begin{equation}\label{eq12}
\|\KK f\|_{p,\varphi;Q}\leq C\|f\|_{p,\varphi;Q},\quad \|\CC[a,f]\|_{p,\varphi;Q}\leq C\|a\|_\ast
\|f\|_{p,\varphi;Q}
\end{equation}
with $C$  independent of $a$ and $f$.
\end{cor}
\proof
Define the extensions
$$
\ol\K(x,\xi)=\begin{cases}
\K(x,\xi)   & (x,\xi)\in Q\times \R^{n+1}_+\sm\{0\}\\
0 & \text{ elsewhere }
\end{cases},\quad
\ol f(x)=\begin{cases}
f(x) &  x\in Q\\
0  &  x\not\in Q.
\end{cases}
$$
Denote by $\ol\KK  f$ the singular integral with a kernel $\ol\K$ and potential $\ol f.$
 Then
$$
|\KK f|\leq |\ol\KK f|\leq C\int_{\R^{n+1}}\frac{|\ol f(y)|}{\rho(x-y)^{n+2}}\,dy
$$
and
$$
\| \KK f\|_{p,\varphi;Q}\leq \|\ol\KK f\|_{p,\varphi;\R^{n+1}}\leq
 C\|\ol f\|_{p,\varphi;\R^{n+1}}= C\|f\|_{p,\varphi;Q}.
$$
The estimate for the commutator follows in a similar way.
\endproof
\begin{cor} \label{locest2}
Let   $a\in VMO$ and  $(p,\varphi)$ be as in Theorem~\ref{th2}.
Then for any $\varepsilon>0$ there exists a positive number
$r_0=r_0(\varepsilon,\eta_a)$ such that for any  $\E_r(x_0)$
with a radius $r\in(0,r_0)$
and all $f\in M_{p,\varphi}(\E_r(x_0))$
\begin{equation}\label{normB}
\|\CC[a,f]\|_{p,\varphi;\E_r(x_0)}
\leq  C\varepsilon\|f\|_{p,\varphi;\E_r(x_0)}
\end{equation}
where $C$ is independent of $\varepsilon$, $f,$  $r$ and $x_0.$
\end{cor}
\proof
Since any $VMO$ function can be   approximated by BUC  functions (see \cite{ChFraL1, Sar}) for each $\varepsilon >0$  there exists $r_0(\varepsilon,\eta_a)$ and $g\in BUC$ with modulus of continuity $\omega_g(r_0)<\varepsilon/2$ such that
$\|a-g\|_\ast<\varepsilon/2. $  Fixing $\E_r(x_0)$
 with $r\in(0,r_0)$ define the function
$$
h(x)=\begin{cases}
\ds g(x) & x\in\E_r(x_0) \\
\ds g\big(x_0+r\frac{x'-x_0'}{\rho(x-x_0)}, t_0+r^2\frac{t-t_0}{\rho^2(x-x_0)}  \big) & x\in \E_r^c(x_0)
\end{cases}
$$
such that $h\in BUC(\R^{n+1})$ and $\omega_h(r_0)\leq \omega_g(r_0)<\varepsilon/2.$     Hence
\begin{align*}
\|\CC[a,f]\|_{p,\varphi;\E_r(x_0)} & \leq \|\CC[a-g,f]\|_{p,\varphi;\E_r(x_0)}+\|\CC[g,f]\|_{p,\varphi;\E_r(x_0)}\\
&  \leq C\|a-g\|_\ast\|f\|_{p,\varphi;\E_r(x_0)} +\|\CC[h,f]\|_{p,\varphi;\E_r(x_0)}
< C\varepsilon \|f\|_{p,\varphi;\E_r(x_0)}\,.
\end{align*}
\endproof

For any $x'\in {\mathbb R}^n_+$ and any fixed $t>0$   define the {\it generalized reflection}
\begin{equation}\label{GR}
\T(x)=(\T'(x),t) \qquad \T'(x) =x'-2x_n\frac{{\bf a}^n(x',t)}{a^{nn}(x',t)}
\end{equation}
where ${\bf a}^n(x)$ is  the last row of the coefficients matrix ${\bf a}(x)$ of \eqref{CDP}.
The function  $\T'(x)$ maps $\R^n_+$ into $\R^n_-$ and the kernel
$\K(x,\T(x)-y)=\K(x, \T'(x)-y',t-\tau)$ is nonsingular  one  for any $x,y\in \D^{n+1}_+.$  Taking $\wtl x\in \D^{n+1}_-$ there  exist positive constants $\kappa_1$ and $\kappa_2$ such that
\begin{equation}\label{CTC}
\kappa_1\rho(\wtl x - y) \leq \rho({\T}(x)-y) \leq \kappa_2 \rho(\wtl x -y).
\end{equation}
For any
  $f\in M_{p,\varphi}({\D}^{n+1}_+)$ and  $a\in BMO({\D}^{n+1}_+)$ define the nonsingular integral operators
\begin{equation}\label{KCf}
\begin{cases}
\ds\wtl   {\KK}	f(x) =\int_{{\D}^{n+1}_+} \K (x,{\T}(x)-y)f(y) dy\\
\ds \wtl   {\CC}  [a,f](x)=\int_{{\D}^{n+1}_+} \K(x,{\T}(x)-y)[a(y)-a(x)]f(y) dy.
\end{cases}
\end{equation}
Since $\K(x, \T(x)-y)$ is still homogeneous  one and satisfies the conditin $b)$ in Definition~\ref{CZK}
we have
$$
|\K(x, \T(x)-y)|\leq \frac{M}{\rho(\T(x)-y)^{n+2}}\leq  \frac{C}{\rho(\wtl x-y)^{n+2}}.
$$
Hence the  operators   \eqref{KCf} are sublinear  and bounded in
 $L_p(\D^{n+1}_+),$ $p\in(1,\infty)$ (cf. \cite{BC}).
The following estimates are simple consequence of  the results in  \S\ref{sec4}.
\begin{thm}\label{nonsing}
Let  $a \in BMO(\D^{n+1}_+)$   and   $f\in M_{p,\varphi}(\D^{n+1}_+) $ with   $(p,\varphi)$ as in Theorem~\ref{th2}. Then the operators  $\wtl\KK f$ and $\wtl\CC[a, f]$ are continuous in $M_{p,\varphi}(\D^{n+1}_+)$  and \begin{equation}\label{KC}
\|\wtl\KK f\|_{p,\varphi;\D^{n+1}_+} \leq   C  \|f\|_{p,\varphi;\D^{n+1}_+},
\quad  \|\wtl\CC[a, f]\|_{p,\varphi;\D^{n+1}_+} \leq   C \|a\|_\ast \,
\|f\|_{p,\varphi;\D^{n+1}_+}
\end{equation}
with a constant independend of $a$ and $f.$
\end{thm}
\begin{cor} \label{localnonsing}
Let  $a\in VMO$ and $(p,\varphi)$ be as above.
Then for any $\varepsilon>0$ there exists a positive number
$r_0=r_0(\varepsilon,\eta_a)$ such that for any  $\E_r^+(x^0)$
with a radius $r\in(0,r_0)$
and all $f\in M_{p,\varphi}(\E_r^+(x^0))$
\begin{equation}\label{tlK}
\|\CC[a,f]\|_{p,\varphi;\E^+_r(x^0)}
\leq  C\varepsilon\|f\|_{p,\varphi;\E^+_r(x^0)},
\end{equation}
where $C$ is independent of
$\varepsilon$, $f,$  $r$ and $x^0$.
\end{cor}

\section{Proof of the main result}\label{sec6}

Consider the problem \eqref{CDP} with $f\in M_{p,\varphi}(Q),$ $(p,\varphi)$ as in Theorem~\ref{th2}.  Since $M_{p,\varphi}(Q)$ is a proper subset of $L_p(Q)$  than \eqref{CDP} is uniquely solvable  and the solution $u$ belongs at least  to $ \overset{\circ}{W}{}^{2,1}_{p}(Q).$  Our aim is to show that this solution  belongs also to $ \overset{\circ}{W}{}^{2,1}_{p,\varphi}(Q).$
For this goal we need a priori estimate of $u$ that we are going to prove in two steps.

{\it  Interior estimate.}
For any  $x_0\in \R^{n+1}_+$ consider  the parabolic semi-cylinders
$\C_r(x_0)=\B_r(x_0')\times (t_0-r^2,t_0).$  Let $v\in C_0^\infty(\C_r)$ and suppose that $v(x,t)=0$ for $t\leq 0.$ According to \cite[Theorem~1.4]{BC}  for any  $x\in \supp\, v$ the following representation formula for the second derivatives of $v$ holds true
\begin{align}
\nonumber
D_{ij}v(x)=& P.V.\int_{\R^{n+1}}\Gamma_{ij}(x,x-y)[a^{hk}(y)-a^{hk}(x)]D_{hk}v(y)dy\\
\label{RF}
&+ P.V. \int_{\R^{n+1}}\Gamma_{ij}(x,x-y)\P v(y)dy+\P v(x)\int_{\SF^n}\Gamma_j(x,y)\nu_id\sigma_y,
\end{align}
where $\nu(\nu_1,\ldots,\nu_{n+1})$ is the outward normal to $\SF^n.$  Here $\Gamma(x,\xi)$ is the fundamental solution of the operator $\P$ and $\Gamma_{ij}(x,\xi)=\partial^2
\Gamma(x,\xi)/\partial\xi_i\partial\xi_j.$   Since any function $v\in W^{2,1}_p$ can be approximated by $C_0^\infty$ functions, the representation formula  \eqref{RF}
still holds for any  $v\in W^{2,1}_{p}(\C_r(x_0)).$
 The properties of the fundamental  solution (cf. \cite{BC,LSU, Sf1}) imply  $\Gamma_{ij}$ are variable Calder\'on-Zygmund kernels in the sense of Definition~\ref{CZK}. Using the notations \eqref{sing} we can write
\begin{align}\nonumber
D_{ij}v(x)=& \CC_{ij}[a^{hk},D_{hk}v](x)\\
\label{RF2}
& + \KK_{ij}(\P v)(x)+\P v(x)\int_{\SF^{n}} \Gamma_j(x,y)\nu_i d\sigma_y\,.
\end{align}
The integrals  $\KK_{ij}$ and $ \CC_{ij}$ are  defined by \eqref{sing} with kernels $\K(x,x-y)=\Gamma_{ij}(x,x-y).$
Because of Corollaries~\ref{locest1} and~\ref{locest2} and the equivalence of the metrics  we get
\begin{equation}\label{eq13}
\|D^2v\|_{p,\varphi;\C_r(x_0)}\leq C(\varepsilon\|D^2v\|_{p,\varphi;\C_r(x_0)}+
\|\P u\|_{p,\varphi;\C_r(x_0)}  )
\end{equation}
for some $r$  small enough. Moving  the norm of $D^2v$ on the left-hand side we get
$$
\|D^2 v\|_{p,\varphi;\C_r(x_0)}\leq C(n,p,\eta_a(r), \|D\Gamma\|_{\infty,Q}  )\|\P v\|_{p,\varphi;\C_r(x_0)}.
$$
Define  a  cut-off function $\phi(x)=\phi_1(x')\phi_2(t),$   with $\phi_1\in C_0^\infty(\B_r(x'_0)),$
$\phi_2\in C_0^\infty(\R) $ such that
$$
\phi_1(x')=\begin{cases}
1 & x'\in \B_{\theta r}(x_0')\\
0 & x'\not\in \B_{\theta'r}(x_0')
\end{cases},\qquad
\phi_2(t)=\begin{cases}
1 & t\in (t_0-(\theta r)^2, t_0]\\
0 & t< t_0-(\theta' r)^2
\end{cases}
$$
with  $\theta\in(0,1),$ $\theta'=\theta(3-\theta)/2>\theta$ and
$|D^s\phi|\leq C [\theta(1-\theta)r]^{-s},$ $s=0,1,2,$ \
$|\phi_t|\sim |D^2\phi|.$  For any solution $u\in W^{2,1}_p(Q) $ of
\eqref{CDP}
define  $v(x)=\phi(x) u(x)\in W^{2,1}_p(\C_r).$ Hence
\begin{align*}
\|D^2u\|_{p,\varphi;\C_{\theta r}(x_0)} & \leq  \| D^2 v\|_{p,\varphi;\C_{\theta'r}(x_0)}\leq
 C\|\P v \|_{p,\varphi;\C_{\theta' r}(x_0)}\\
& \leq  C\left(\|f\|_{p,\varphi;\C_{\theta' r}(x_0)}+\frac{\|Du\|_{p,\varphi;\C_{\theta' r}(x_0)}}{\theta(1-\theta)r}
 +\frac{\|u\|_{p,\varphi;\C_{\theta' r}(x_0)}}{[\theta(1-\theta)r]^2}  \right).
\end{align*}
Hence
\begin{align*}
\big[\theta(1-\theta)r\big]^2 &\| D^2u\|_{p,\varphi;\C_{\theta r}(x_0)}\\
\leq &\  \left([\theta(1-\theta)r]^2\|f\|_{p,\varphi;\C_{\theta' r}(x_0)}+\theta(1-\theta)r
\|Du \|_{p,\varphi;\C_{\theta' r}(x_0)}+\|u\|_{p,\varphi;\C_{\theta' r}(x_0)}  \right)\\
&\  \left(\text{by the choice   of  } \theta' \text{ it follows }  \theta(1-\theta)\leq 2 \theta'(1-\theta')\right)\\
\leq &\  C \left( r^2\|f\|_{p,\varphi;Q}+\theta'(1-\theta')r\|Du \|_{p,\varphi;\C_{\theta' r}(x_0)}+
\|u\|_{p,\varphi;\C_{\theta' r}(x_0)}    \right)\,.
\end{align*}

Introducing the semi-norms
$$
\Theta_s=\sup_{0<\theta<1} \big[\theta(1-\theta)r \big]^s \|D^s u \|_{p,\varphi;\C_{\theta r}(x_0)}\qquad s=0,1,2
$$
the above inequality becomes
\begin{equation}\label{theta}
[\theta(1-\theta)r]^2\|D^2u\|_{p,\varphi;\C_{\theta r}(x_0)} \leq \Theta_2
 \leq C
\left(r^2\|f\|_{p,\varphi;Q} +\Theta_1+\Theta_0 \right)\,.
\end{equation}
The  interpolation inequality  \cite[Lemma~4.2]{Sf3} gives that there  exists a positive constant $C$  independent of $r$ such that
$$
\Theta_1\leq \varepsilon\, \Theta_2+\frac{C}{\varepsilon}\, \Theta_0\qquad \text{ for any } \varepsilon\in(0,2).
$$
Thus   \eqref{theta} becomes
$$
[\theta(1-\theta)r]^2\|D^2u\|_{p,\varphi;\C_{\theta r}(x_0)}\leq\Theta_2\leq
 C\left(r^2\|f\|_{p,\varphi;Q}+\Theta_0 \right) \quad \forall\  \theta\in(0,1).
$$
Taking  $\theta =1/2$ we get  the    Caccioppoli-type estimate
$$
\|D^2u\|_{p,\varphi; \C_{r/2}(x_0)} \leq  C\left( \|f\|_{p,\varphi;Q}+
 \frac1{r^2}\|u\|_{p,\varphi;\C_r(x_0)} \right).
 $$
To estimate  $u_t $ we exploit the parabolic structure of  the equation
and the boundedness of the coefficients
\begin{align*}
\|u_t\|_{p,\varphi;\C_{r/2}(x_0)}& \leq \|{\bf a}\|_{\infty;Q} \|D^2u\|_{p,\varphi;\C_{r/2}(x_0)}
+\|f\|_{p,\varphi;\C_{r/2}(x_0)}\\
&  \leq C\big( \|f\|_{p,\varphi;Q} + \frac1{r^2}\|u\|_{p,\varphi;\C_r(x_0)} \big).
\end{align*}
 Consider  cylinders $Q'=\Omega'\times(0,T)$ and $Q''=\Omega''\times(0,T)$ with $\Omega'\subset\subset\Omega''\subset\subset\Omega,$
by standard covering procedure and partition of the unity we get
\begin{equation}\label{intest}
\|u\|_{W^{2,1}_{p,\varphi}(Q')}  \leq C
     \big(\|f\|_{p,\varphi;Q}
+\|u\|_{p,\varphi;Q''}\big)
\end{equation}
where  $C$ depends on $n, p, \Lambda, T,
\|D\Gamma\|_{\infty;Q}, \eta_{\bf a}(r),$ $\|{\bf a}\|_{\infty,Q}$ and $\text{dist}(\Omega',\partial\Omega'').$

{\it Boundary estimates.}
For any fixed $(x^0,r)\in \R^{n-1}\times\R_+$ define the semi-cylinders
$$
\C_r^+(x^0)=\B_r^+({x^0}')\times (0,r^2)=\{ |x^0-x'|<r, x_n>0, 0<t<r^2 \}
$$
with $\SS_r^+=\{(x'',0,t): \    |x^0-x''|<r, 0<t<r^2 \}.$ For any solution $u\in W^{2,1}_p(\C_r^+(x^0))$ with
$\supp u\in \C_r^+(x^0)$ the following boundary representation formula holds (cf. \cite{BC})
\begin{align*}
D_{ij}u(x)=& \CC_{ij}[a^{hk},D_{hk}u](x) +\KK_{ij} (\P u)(x)\\
& + \P u(x) \int_{\SF^{n}}\Gamma_j(x,y)\nu_i d\sigma_y - \II_{ij}(x)
\end{align*}
where
\begin{align*}
&\II_{ij}(x)=\wtl\KK_{ij}(\P u)(x)+ \wtl\CC_{ij}[a^{hk},D_{hk}u](x), \quad i,j=1,\ldots, n-1,\\
&\II_{in}(x)=\II_{ni}(x)=\sum_{l=1}^n\left( \frac{\partial\T(x)}{\partial x_n} \right)^l \left[\wtl\CC_{il}[a^{hk},D_{hk}u](x)+\wtl\KK_{il}(\P u)(x) \right], \quad  i=1,\ldots,n-1,\\
&\II_{nn}(x)=\sum_{r,l=1}^n \left( \frac{\partial\T(x)}{\partial x_n} \right)^r \left( \frac{\partial\T(x)}{\partial x_n} \right)^l   \left[\wtl\CC_{rl}[a^{hk},D_{hk}u](x)+\wtl\KK_{rl}(\P u)(x) \right],\\
&\frac{\partial \T(x)}{\partial x_n}=\left(-2\frac{a^{n1}(x)}{a^{nn}(x)},\ldots,
-2\frac{a^{n n-1}(x)}{a^{nn}(x)},-1,0    \right)\,.
\end{align*}
Here $\wtl\KK_{ij}$ and $ \wtl\CC_{ij}$   are the operators defined by \eqref{KCf} with kernels $\K(x,\T(x)-y)=\Gamma_{ij}(x,\T(x)-y).$
Applying the estimates \eqref{KC} and \eqref{tlK} and having in mind that the components of the vector $\frac{\partial\T(x)}{\partial x_n}$ are bounded  we get
$$
\|D^2u\|_{p,\varphi;\C^+_r(x^0)}\leq C\big( \|\P u\|_{p,\varphi;\C_r^+(x^0)}+\|u\|_{p,\varphi;\C_r^+(x^0)}  \big).
$$
The Jensen inequality applied to $u(x)=\int_0^t u_s(x',s)ds$  and the parabolic structure of the equation give
$$
\|u\|_{p,\varphi;\C_r^+(x^0)}\leq C r^2\|u_t\|_{p,\varphi;\C_r^+(x^0)}\leq C(\|f\|_{p,\varphi;Q} +
r^2\|u\|_{p,\varphi;\C_r^+(x^0)}  ).
$$
Taking $r $ small enough we can move the norm of $u$ on the left-hand side obtaining
$$
\|u\|_{p,\varphi;\C_r^+}\leq C\|f\|_{p,\varphi;Q}
$$
with a constant $C$ depending on $n,p,\Lambda, T, \eta_{\ba}, \|\ba\|_{\infty,Q}.$ 
By covering of the boundary with small cylinders, partition of the unit subordinated of that covering and local flattering of $\partial\Omega$ we get that
\begin{equation}\label{boundest}
\|u\|_{W^{2,1}_{p,\varphi}(Q\sm Q')}\leq C\|f\|_{p,\varphi;Q}\,.
\end{equation}
Unifying \eqref{intest} and \eqref{boundest} we get \eqref{apriori}.


\medskip

\begin{thebibliography}{99}





\bibitem{A}P.  Acquistapace,  {\em  On $BMO$ regularity for linear elliptic systems},
Ann. Mat. Pura Appl.,  {\bf 161},  231--270, 1992.

\bibitem{AGM}A.  Akbulut, V.S.  Guliyev, R.  Mustafayev,  {\em 
On the boundedness of the maximal operator and singular
integral operators in generalized Morrey spaces},  Math. Bohem., {\bf 137} (1),  27--43, 2012.

\bibitem{BC}M. Bramanti,  M.C.   Cerutti,   {\em $W_{p}^{1,2}$ solvability
for the Cauchy--Dirichlet problem for parabolic equations with VMO
 coefficients}, Comm.  Partial Diff. Eq., {\bf 18}, 1735--1763, 1993.


\bibitem{CarPickSorStep}M.  Carro,  L. Pick, J.  Soria, V.D.  Stepanov,  {\em 
 On embeddings between classical Lorentz spaces},
Math. Inequal. Appl., {\bf  4} (3),  397--428, 2001.

\bibitem{ChFra}F. Chiarenza, M.  Frasca, 
{\em  Morrey spaces and Hardy-Littlewood maximal function},
Rend. Mat., {\bf  7}, 273--279, 1987.

\bibitem{ChFraL1}F.  Chiarenza, M.   Frasca, M., P.  Longo,  {\em
 Interior $W^{2,p}$-estimates for
nondivergence elliptic equations with discontinuous coefficients},
Ricerche Mat., {\bf  40}, 149--168, 1991.

\bibitem{FR}E.B. Fabes, N.  Rivi\`ere,  {\em Singular integrals with mixed
homogeneity}, Studia Math., {\bf 27},  19--38, 1996.

\bibitem{G1}V.S.  Guliyev,  {\em  Integral operators on function spaces on the homogeneous groups and
  on domains in $\R^n$}. Doctor's degree dissertation,
  Mat. Inst. Steklov, Moscow, 1994, 329 pp. (in Russian)

\bibitem{G2}V.S.  Guliyev, 
{\em  Boundedness of the maximal, potential and singular operators in the generalized Morrey spaces},
J. Inequal. Appl.,  2009, Art. ID 503948, pp. 20.

\bibitem{GM}V.S.   Guliyev, F.M.  Mushtagov, 
 {\em Parabolic equations with $VMO$ coefficients in weighted Lebesgue spaces,}
Proc. Razmadze Math. Inst., {\bf 137}, 1--27, 2005.

\bibitem{GAKS}V.S.  Guliyev, S.S.   Aliyev, T.  Karaman, P.  Shukurov, 
 {\em Boundedness of sublinear operators and commutators on generalized Morrey spaces,}
Integral  Equ. Oper. Theory, {\bf 71} (3),  327--355, 2011.

\bibitem{GS}V.S.  Guliyev, L.G.  Softova,
{\em Global regularity in generalized Morrey spaces of solutions to non-divergence
 elliptic equations with $VMO$ coefficients}, Potential Anal.,  (on line first), DOI 10.1007/s11118-012-9299-4.

\bibitem{JN}   F.  John, L.  Nirenberg, 
  {\em  On functions of bounded mean oscillation},
  Commun. Pure Appl. Math., {\bf 14}, 415--426, 1961.



\bibitem{Jones} P.W.  Jones, 
{\em  Extension theorems for BMO},
Indiana Univ. Math. J., {\bf  29}, 41--66, 1980.

\bibitem{LSU}O.A.  Ladyzhenskaya, V.A.    Solonnikov, N.N.   Ural'tseva, 
  {\em   Linear and Quasilinear Equations of Parabolic Type},  Transl. Math.
    Monographs  23, Amer. Math. Soc., Providence, R.I., 1968.

\bibitem{Mi}T.   Mizuhara,  
 {\em  Boundedness of some classical operators on generalized Morrey spaces}, Harmonic Anal., Proc. Conf., Sendai/Jap. 1990, ICM-90 Satell. Conf. Proc.,  183--189, 1991. 


\bibitem{Mo}C.B.   Morrey,  
{\em   On the solutions of quasi-linear elliptic partial
  differential equations}, Trans. Amer. Math. Soc.,  {\bf 43}, 126--166, 1938.

\bibitem{Na}E.   Nakai,   
 {\em Hardy-Littlewood maximal operator, singular integral
operators and the Riesz potentials on generalized Morrey spaces},   Math. Nachr.,
{\bf 166}, 95--103, 1994.

\bibitem{PRS}D.K.    Palagachev, M.A.  Ragusa, L.G.   Softova, 
{\em Cauchy-Dirichlet problem in Morrey spaces for parabolic equations with discontinuous coefficients}, Bolletino U.M.I., {\bf 8} 6-B, 667--683, 2003.


\bibitem{PS1}D.K.  Palagachev, L.G.   Softova,  
{\em  Singular integral operators, Morrey spaces and fine regularity of solutions to PDE's},
      Potential Anal.,  {\bf 20}, 237--263, 2004. 


\bibitem{Sar}D.  Sarason,  
{\em   On functions of vanishes mean oscillation},
Trans. Amer. Math. Soc., {\bf  207},   391--405, 1975.



\bibitem{Sf2}L.G.  Softova,  
 {\em  Singular integrals and commutators in generalized Morrey spaces}, Acta Math. Sin., Engl. Ser., {\bf 22}, 757--766, 2006.

\bibitem{Sf1}L.G.  Softova, 
 {\em  Singular integral operators in Morrey spaces and interior
regularity of solutions to systems of linear PDE’s}, J. Glob. Optim., {\bf 40}, 427--442, 2008.

\bibitem{Sf3}L.G.   Softova, 
 {\em Morrey-type regularity of solutions to parabolic problems
with discontinuous data}, Manuscr. Math.,  {\bf 136} (3--4),  365--382, 2011.

\bibitem{Sf4}L.G.  Softova,  
{\em The Dirichlet problem for elliptic equations with $VMO$ coefficients in generalized Morrey spaces}, in: Advances in Harmonic Analysis and Operator Theory, The Stefan Samko Anniversary Volume,  (in print), 2012.

\end{thebibliography}
\end{document}